\documentclass[12pt]{article}
\usepackage[namelimits, reqno]{amsmath}
\usepackage{amssymb}
\usepackage{graphicx}
\usepackage{xcolor,graphicx}

\def\R{I\!\!R}

\def\lim{\mathop {\rm lim}}

\begin{document}

\begin{center}\bf
\large WAVE EQUATION WITH THREE-INVERSE  
SQUARE POTENTIAL ON $R_+^3$ 
\end{center}
\begin{center} Yehdhih Mohamed Abdelhaye, Badahi  Mohamed 
and Mohamed Vall Ould Moustapha
\end{center}

\newtheorem{Theorem}{\quad Theorem}[section] 

\newtheorem{Definition}[Theorem]{\quad Definition} 

\newtheorem{Corollary}[Theorem]{\quad Corollary} 

\newtheorem{Lemma}[Theorem]{\quad Lemma} 

\newtheorem{Example}[Theorem]{\quad Example} 
\scriptsize{{\bf Abstract} In this note we give explicit solutions to the wave equation associated to the Schr\"odinger operator  with three-inverse square potential on $R_+^3$.}\\
{\bf Key words }: Three inverse square potential, Cauchy problem, Wave equation, Lauricella hypergeometric functions.

\section{Introduction and statement of results}
The wave equation, the heat equation and the Laplace equation are known as three fundamental equations in partial differential equations and occur in many branches of physics, in applied mathematics and in engineering. In this note we give explicit formulas for the solutions of the following Cauchy problem for the wave equation with three-inverse square potential 
$$(W)_{(\nu, \nu', \nu'')} \left \{\begin{array}{cc}
(a)\ \ \ [ \Delta+v_{\nu, \nu', \nu''}(x)]u(t, p)=\frac{\partial^2}{\partial t^2} u(t, p)&(t, p)\in R\times R_+^{3}\\
(b)\qquad\ \  u(0, p)= 0& \frac{\partial}{ \partial t} u(0, p)= u_1(p), u_1\in
C_0^\infty(R_+^{3})\end{array}
\right.$$
with
$\Delta=\frac{\partial^{2}}{\partial x^{2}}+\frac{\partial^{2}}{\partial y^{2}}+\frac{\partial^{2}}{\partial z^{2}}$
is the  Laplacien of $R^{3}$, the three inverse square potential is given by
$v_{\nu,\nu', \nu''}(x, y, z)= \frac{1/4-\nu^2}{x^2}+\frac{1/4-\nu'^2}{y^2}+\frac{1/4-\nu''^2}{z^2}$
where $\nu, \nu', \nu'' $  are real parameters.
The Cauchy problem for the wave equation with the inverse square potential in
Euclidean space $\mathbb{R}^n$ is extensively studied (Cheeger and Taylor $[3]$)). The bi-inverse square potential has been considered by (Boyer $[2]$)and Ould Moustapha $[7]$.
 The case considered most
frequentely is obviously the one where $(\nu, \nu', \nu'')=(\pm 1/2, \pm 1/2, \pm 1/2)$, the equation in $(W)_{\nu, \nu', \nu''}$
 then turns into the classical wave equation on the Euclidean space
$\mathbb{R}^3$ and this equation appears in several branches of mathematics
and physics (Folland $[5], p.171$).
Now we state the main results of this paper:\\

{Theorem A}{ \it
For $ (t,p,p')\in R_+\times {\mathbb{R}^\ast}^3_+\times{\mathbb{R}^\ast}^3_+$ the
functions:\\
$ W_{(b,b',b'')}(t,p,p')
 =\frac{c_3(xx')^{b}(yy')^{b'}(zz')^{b''}}{\left(t^2-|p+p'|^2\right)^{1+b+b'+b''}}
 F_A^{(3)}\left(1+b+b'+b'', b, b', b'', 2 b, 2 b', 2b'', \frac{-4xx'}{t^2-|p+p'|^2},\frac{-4yy'}{t^2-|p+p'|^2},\frac{-4zz'}{t^2-|p+p'|^2}\right)$\\
with $b \in \{\beta,1-\beta\}$, $b' \in \{\beta',1-\beta'\}$,  $b'' \in \{\beta'',1-\beta''\}$
are independent solutions of the wave equation with three-inverse square
potential on $\R_+^3$ $(a)$ where $\beta=1/2+\nu$ and $\beta'=1/2+\nu'$, $\beta''=1/2+\nu''$ 
and $F_A^{(3)}(a, b, b', b'', c, c', c''; w, w', w'')$ is the three
variables triple series  $F_A^{(3)}$ Lauricella hypergeometric function given by $[1], p.114$
$$
F_A^{(3)}\left(a, b, b', b'', c, c', c'', w, w', w''\right)=\sum_{m, n, p\geq
  0}\frac{(a)_{m+n+p}(b)_m(b')_n(b'')_p} {(c)_m m!(c')_n n!(c'')_p
  p!}w^m w'^n w''^p\eqno(1.1)
$$
{\bf Theorem B}{\it
The Cauchy problem for the wave equation with
three-inverse square potential on the $\R_+^3$ has the
solutions given by:
$$
  u(t, p)=\frac{\partial}{t \partial t}\int_{|p+p'|<t}W_{(b, b' ,b'')}(t, p, p')f(p')dp'
\eqno(1.2)$$
where the kernel $W_{(b, b', b'')}$ is as in the Theorem A and the constant $c_3$ is given by
$$
  c_3=\frac{(-1)^{b+b'+b''} 
  \Gamma(1/2+b+b'+b'')}{\sqrt{\pi}\Gamma(1/2+b)\Gamma(1/2+b')\Gamma(1/2+b'')}
\eqno(1.3)$$
with  $dp'=dx'dy'dz'$ is the Lebesgue measure on
$\mathbb{R}^3$ }\\

\section{Wave equation with three-inverse square potential on Euclidean space
  $\mathbb{R}^3$}
{ \bf Proof of Theorem A}\\
In what follows we give a direct proof of the theorem A.\\
 Let $a=t^2-|p+p'|^2$, $t\in \mathbb{R}$, $p,p'\in {\mathbb{R}^\ast}^3$  set:\\
$\Omega\varphi(t,p)=\left(x x'\right)^{-\beta}\left(y
    y'\right)^{-\beta'}\left(z
    z'\right)^{-\beta''} a^{-\alpha}\times$
  $$\left[\Delta-\frac{\beta(\beta-1)}{x^2}-\frac{\beta'(\beta'-1)}{y^2}-\frac{\beta''(\beta''-1)}{z^2}-\frac{\partial^2}{\partial t^2}\right]\left(xx'\right)^{\beta}\left(yy '\right)^{\beta'}\left(z z '\right)^{\beta''}a^{\alpha}\varphi(t,x)\eqno(2.1)$$
then  we have\\
$ \Omega\varphi(t,p)=
\{\Delta-\frac{\partial^2}{\partial t^2}+[\frac{2\beta}{x}-\frac{4\alpha(x+x')}{a}]\frac{\partial}{\partial x}+[\frac{2\beta'}{y}-\frac{4\alpha(
y+y')}{a}]\frac{\partial}{\partial y}+[\frac{2\beta''}{z}-\frac{4\alpha(
z+z')}{a}]\frac{\partial}{\partial z}-\frac{4\alpha
t}{a}\frac{\partial}{\partial t}$
$$-\frac{4\alpha }{a}[\frac{\beta x'}{x}+
\frac{\beta' y'}{y}+\frac{\beta'' z'}{z}]-\frac{4\alpha
}{a}[\alpha+1+\beta+\beta'+\beta'']\}\varphi(t, x)\eqno(2.2)$$
Now set 
$$
w=\frac{-4x x'}{t^2-|p+p'|^2}, w'=\frac{-4y y'}{t^2-|p+p'|^2}, w''=\frac{-4z z'}{t^2-|p+p'|^2}\eqno(2.3)$$
we can write:\\ $$
  \frac{\partial}{\partial x}=\frac{\partial w}{\partial x}\frac{\partial}{\partial w}+\frac{\partial w'}{\partial x}\frac{\partial}{\partial w'}+\frac{\partial w''}{\partial x}\frac{\partial}{\partial w''}
;
  \frac{\partial}{\partial y}=\frac{\partial w}{\partial y}\frac{\partial}{\partial w}+\frac{\partial w'}{\partial y}\frac{\partial}{\partial w'}+\frac{\partial w''}{\partial y}\frac{\partial}{\partial w''}\eqno(2.4)$$
 $$ \frac{\partial}{\partial z}=\frac{\partial w}{\partial z}\frac{\partial}{\partial w}+\frac{\partial w'}{\partial z}\frac{\partial}{\partial w'}+\frac{\partial w''}{\partial z}\frac{\partial}{\partial w''};
 \frac{\partial}{\partial t}=\frac{\partial
  w}{\partial t}\frac{\partial}{\partial w}+\frac{\partial w'}{\partial
  t}\frac{\partial}{\partial w'}+\frac{\partial w''}{\partial
  t}\frac{\partial}{\partial w''}\eqno(2.5)$$
We have:\\
$
  \Omega\varphi(t,p)=
  \left[\left(\frac{\partial w}{\partial x}\right)^2+\left(\frac{\partial w}{\partial y}\right)^2
    +\left(\frac{\partial w}{\partial z}\right)^2
    -\left(\frac{\partial w}{\partial t}\right)^2\right]
  \frac{\partial^2}{\partial w^2}
  \\+\left[\left(\frac{\partial w'}{\partial x}\right)^2
    +\left(\frac{\partial w'}{\partial y}\right)^2+\left(\frac{\partial w'}{\partial z}\right)^2
    -\left(\frac{\partial w'}{\partial t}\right)^2\right]
  \frac{\partial^2}{\partial w'^2}+\left[\left(\frac{\partial w''}{\partial x}\right)^2
    +\left(\frac{\partial w''}{\partial y}\right)^2+\left(\frac{\partial w''}{\partial z}\right)^2
    -\left(\frac{\partial w''}{\partial t}\right)^2\right]
  \frac{\partial^2}{\partial w''^2}
  \\+2\left[\frac{\partial w}{\partial x}\frac{\partial w'}{\partial x}
    +\frac{\partial w}{\partial y}\frac{\partial w'}{\partial y}+\frac{\partial w}{\partial z}\frac{\partial w'}{\partial z}-\frac{\partial w}{\partial t}\frac{\partial w'}{\partial t}\right]\frac{\partial^2}{\partial w\partial w'}+2\left[\frac{\partial w}{\partial x}\frac{\partial w''}{\partial x}
    +\frac{\partial w}{\partial y}\frac{\partial w''}{\partial y}+\frac{\partial w}{\partial z}\frac{\partial w''}{\partial z}-\frac{\partial w}{\partial t}\frac{\partial w''}{\partial t}\right]\frac{\partial^2}{\partial w \partial w''}\\+2\left[\frac{\partial w'}{\partial x}\frac{\partial w''}{\partial x}
    +\frac{\partial w'}{\partial y}\frac{\partial w''}{\partial y}+\frac{\partial w'}{\partial z}\frac{\partial w''}{\partial z}-\frac{\partial w'}{\partial t}\frac{\partial w''}{\partial t}\right]\frac{\partial^2}{\partial w'\partial w''}
    +\left[\frac{\partial^2 w}{\partial x^2}
    +\frac{\partial^2 w}{\partial y^2}+\frac{\partial^2 w}{\partial z^2}
    -\frac{\partial^2 w}{\partial t^2}\right]
  \frac{\partial}{\partial w}
  \\+[\frac{\partial^2 w'}{\partial x^2}+\frac{\partial^2 w'}{\partial y^2}
  +\frac{\partial^2 w'}{\partial z^2}
  -\frac{\partial^2 w'}{\partial t^2}]\frac{\partial}{\partial w'}
+[\frac{\partial^2 w''}{\partial x^2}+\frac{\partial^2 w''}{\partial y^2}
  +\frac{\partial^2 w''}{\partial z^2}
  -\frac{\partial^2 w''}{\partial t^2}]\frac{\partial}{\partial w''}
  \\+[A_x \frac{\partial w}{\partial x}+A_y \frac{\partial w}{\partial
    y}+A_z \frac{\partial w}{\partial
    z}-A_t\frac{\partial w}{\partial t}]\frac{\partial }{\partial w}+[A_x
  \frac{\partial w'}{\partial x}+A_y \frac{\partial w'}{\partial
    y}+A_z \frac{\partial w'}{\partial
    z}-A_t\frac{\partial w'}{\partial t}]\frac{\partial }{\partial
    w'}\\
+[A_x
  \frac{\partial w''}{\partial x}+A_y \frac{\partial w''}{\partial
    y}+A_z \frac{\partial w''}{\partial
    z}-A_t\frac{\partial w''}{\partial t}]\frac{\partial }{\partial
    w''}+\frac{4\alpha}{a}[\frac{\beta x'}{x}+\frac{\beta'
    y'}{y}+\frac{\beta''
    z'}{z}] +\frac{4\alpha}{a}\left[\alpha+1+\beta+\beta'+\beta''\right]
\varphi(z, z')\ \ \ \ \ \  \ \  \ \ \  \ \ \  \ \  (2.6)
$\\
where

   $$A_x=\frac{2\beta}{x}-\frac{4\alpha(x+x')}{a},
   A_y=\frac{2\beta'}{y}-\frac{4\alpha(y+y')}{a}, A_z=\frac{2\beta''}{z}-\frac{4\alpha(z+z')}{a},A_t=\frac{4\alpha t}{a}\eqno(2.7)
$$
 We have:
$$
 \frac{\partial w}{\partial x}=\frac{-4x' a-8(x+x')x x'}{a^2}; \frac{\partial w'}{\partial x}=\frac{-8y y'(x+x')}{a^2};\frac{\partial w''}{\partial x}=\frac{-8z z'(x+x')}{a^2}
 \eqno(2.8)$$
$$
 \frac{\partial w}{\partial y}=\frac{-8x x'(y+y')}{a^2}; \frac{\partial w'}{\partial y}=\frac{-4y'a-8y y'(y+y')}{a^2};\frac{\partial w''}{\partial y}=\frac{-8z z'(y+y')}{a^2}
 \eqno(2.9)$$
$$
 \frac{\partial w}{\partial z}=\frac{-8x x'(z+z')}{a^2}; \frac{\partial w'}{\partial z}=\frac{-8y y'(z+z')}{a^2};\frac{\partial w''}{\partial z}=\frac{-4z'a-8z z'(z+z')}{a^2}
 \eqno(2.10)$$
$$
  \frac{\partial w}{\partial t}=\frac{8x x' t}{a^2}; \frac{\partial w'}{\partial t}=\frac{8y y' t}{a^2}, \frac{\partial w''}{\partial t}=\frac{8z z' t}{a^2}
\eqno(2.11)$$
  
$ \frac{\partial^2 w}{\partial x^2}=\frac{-24x x' a^2-16x'^2a^2-32(x+x')^2x x'
    a}{a^4};
  \frac{\partial^2 w}{\partial y^2}= \frac{-8y y' a^2-32x x'(y+y')^2a}{a^4};
  \frac{\partial^2 w}{\partial z^2}= \frac{-8z z' a^2-32x x'(z+z')^2a}{a^4}\ \ \ \ \ \ \ (2.12)
$\\

  $
  \frac{\partial^2 w'}{\partial x^2}= \frac{-8x x' a^2-32 y y'(x+x')^2a}{a^4}
;\frac{\partial^2 w'}{\partial y^2}=\frac{-24y y' a^2-16y'^2a^2-32(y+y')^2y y' a}{a^4}
; \frac{\partial^2 w'}{\partial z^2}=\frac{-8y y' a^2-32 a zz'(y+y')^2a}{a^4}
(2.13)$\\

$
\frac{\partial^2 w''}{\partial x^2}= \frac{-8x x' a^2-32 z z'(x+x')^2a}{a^4}
;\frac{\partial^2 w''}{\partial y^2}=\frac{-8y y' a^2-32zz'(y+y')^2a}{a^4}
 ; \frac{\partial^2 w''}{\partial z^2}=\frac{-24z z' a^2-16z'^2a^2-32(z+z')^2z z' a}{a^4}
(2.14)$

$$
  \frac{\partial^2 w}{\partial t^2}=\frac{8x x' a^2-32a x x't^2}{a^4};
  \frac{\partial^2 w'}{\partial t^2}=\frac{8y y' a^2-32a y y't^2}{a^4},
  \frac{\partial^2 w''}{\partial t^2}=\frac{8z z' a^2-32a z z't^2}{a^4}\eqno(2.15)$$

from $(2.8)-(2.11)$ we have
$$
  \left(\frac{\partial w}{\partial x}\right)^2+\left(\frac{\partial w}{\partial y}\right)^2+\left(\frac{\partial w}{\partial z}\right)^2-\left(\frac{\partial w}{\partial t}\right)^2=\frac{w^2}{x^2}(1-w)\eqno(2.16)
$$
$$
  \left(\frac{\partial w'}{\partial x}\right)^2+\left(\frac{\partial
      w'}{\partial y}\right)^2+\left(\frac{\partial w'}{\partial z}\right)^2-\left(\frac{\partial w'}{\partial
      t}\right)^2=\frac{w'^2}{y^2}(1-w')
\eqno(2.17)$$
$$
  \left(\frac{\partial w''}{\partial x}\right)^2+\left(\frac{\partial
      w''}{\partial y}\right)^2+\left(\frac{\partial w''}{\partial z}\right)^2-\left(\frac{\partial w''}{\partial
      t}\right)^2=\frac{w''^2}{z^2}(1-w'')
\eqno(2.18)$$

$$
  2\left[\frac{\partial w}{\partial x}\frac{\partial w'}{\partial x}
    +\frac{\partial w}{\partial y}\frac{\partial w'}{\partial y}+\frac{\partial w}{\partial z}\frac{\partial w'}{\partial z}    -\frac{\partial w}{\partial t}\frac{\partial w'}{\partial t}\right]
  =\frac{w^2}{x^2}w'+\frac{w'^2}{y^2}w
\eqno(2.19)$$
$$
  2\left[\frac{\partial w}{\partial x}\frac{\partial w''}{\partial x}
    +\frac{\partial w}{\partial y}\frac{\partial w''}{\partial y}+\frac{\partial w}{\partial z}\frac{\partial w''}{\partial z}    -\frac{\partial w}{\partial t}\frac{\partial w''}{\partial t}\right]
  =\frac{w^2}{x^2}w''+\frac{w''^2}{z^2}w
\eqno(2.20)$$
$$
  2\left[\frac{\partial w'}{\partial x}\frac{\partial w''}{\partial x}
    +\frac{\partial w'}{\partial y}\frac{\partial w''}{\partial y}+\frac{\partial w'}{\partial z}\frac{\partial w''}{\partial z}    -\frac{\partial w'}{\partial t}\frac{\partial w''}{\partial t}\right]
  =\frac{w'^2}{y^2}w''+\frac{w''^2}{z^2}w'
\eqno(2.21)$$
from $(2.12)-(2.15)$ we have
$$
  \frac{\partial^2 w}{\partial x^2}+ \frac{\partial^2 w}{\partial
    y^2}+\frac{\partial^2 w}{\partial z^2}-\frac{\partial^2 w}{\partial t^2}=\frac{-w^2}{x^2}+\frac{2}{a}w
\eqno(2.22)$$
$$
  \frac{\partial^2 w'}{\partial x^2}+ \frac{\partial^2 w'}{\partial
    y^2}+\frac{\partial^2 w'}{\partial z^2}-\frac{\partial^2 w'}{\partial t^2}=\frac{-w'^2}{y^2}+\frac{2}{a}w'
\eqno(2.23)$$
$$
  \frac{\partial^2 w''}{\partial x^2}+ \frac{\partial^2 w''}{\partial
    y^2}+\frac{\partial^2 w''}{\partial z^2}-\frac{\partial^2 w''}{\partial t^2}=\frac{-w''^2}{z^2}+\frac{2}{a}w''
\eqno(2.24)$$
from $(2.7)-(2.11)$ we have
$$
  A_x\frac{\partial w}{\partial x}+A_y\frac{\partial w}{\partial y}+A_z\frac{\partial w}{\partial z}-A_t\frac{\partial w}{\partial t}=
\frac{4\alpha+4(\beta+\beta'+\beta'')}{a}w+2\beta\frac{w}{x^2}+\alpha\frac{w^2}{x^2}-\beta\frac{w^2}{x^2}+(\frac{
  \beta' w'}{y^2}+\frac{\beta'' w''}{z^2})w
\eqno(2.25)$$
$$
  A_x\frac{\partial w'}{\partial x}+A_y\frac{\partial w'}{\partial y}+A_z\frac{\partial w'}{\partial z}-A_t\frac{\partial w'}{\partial t}=
\frac{4\alpha+4(\beta+\beta'+\beta'')}{a}w'+2\beta\frac{w'}{y^2}+\alpha\frac{w'^2}{y^2}-\beta\frac{w'^2}{y^2}+(\frac{
  \beta w}{x^2}+\frac{\beta'' w''}{z^2})w'
\eqno(2.26)$$

$$
  A_x\frac{\partial w''}{\partial x}+A_y\frac{\partial w''}{\partial y}+A_z\frac{\partial w''}{\partial z}-A_t\frac{\partial w''}{\partial t}=
\frac{4\alpha+4(\beta+\beta'+\beta'')}{a}w''+2\beta''\frac{w''}{z^2}+\alpha\frac{w''^2}{z^2}-\beta''\frac{w''^2}{z^2}+(\frac{
  \beta w}{x^2}+\frac{\beta' w'}{y^2})w''
\eqno(2.27)$$
To replace in the formula $(2.6)$ using the
formulas $(2.16)-(2.27)$ we get:\\
$
 \Omega\varphi=wx^{-2}A_{\alpha,\beta}(w, w',w'')\varphi+w'y^{-2}A_{\alpha,\beta'}(w', w, w'')\varphi +w''z^{-2}A_{\alpha,\beta''}(w'', w, w')\varphi +$$$
\frac{4}{a}(\alpha+\beta+\beta'+\beta''+1)\left[\left(
  w\frac{\partial }{\partial w}+w'\frac{\partial }{\partial
    w'}+ w''\frac{\partial }{\partial w''}\right)\varphi(w,w',w'')-\varphi(w, w',w'')\right]\eqno(2.28)$$
Take $\alpha=-1-\beta-\beta'-\beta''$ we get
$ \Omega\varphi=0$  is equivalent to
$$
wx^{-2}A_{\alpha,\beta}(w, w',w'')\varphi+w'y^{-2}A_{\alpha,\beta'}(w', w, w'')\varphi +w''z^{-2}A_{\alpha,\beta''}(w'', w, w')\varphi=0 \eqno(2.29)
$$
with\\
$
A_{\alpha,\beta}(w, w', w'')\varphi(w, w', w'')=
   [w(1-w)\frac{\partial^2}{\partial w^2}-w( w'\frac{\partial^2}{\partial
     w' \partial
     w}+w''\frac{\partial^2}{\partial
     w'' \partial
     w})+$ $$ +\left[2\beta+(-\alpha+\beta+1)w\right]\frac{\partial}{\partial
     w}-\beta' w' \frac{\partial}{\partial w'}-\beta'' w'' \frac{\partial}{\partial w''}+\alpha\beta]\varphi(w, w', w'')\eqno(2.30)$$
From the formula $(2.29)$ we have
$$A_{\alpha,\beta}(w, w', w'')\varphi(w, w', w'')=A_{\alpha,\beta'}(w', w, w'')\varphi(w, w', w'')=A_{\alpha,\beta''}(w, w', w)\varphi(w, w', w'')=0\eqno(2.31)$$

$ \begin{aligned}
   &[w(1-w)\frac{\partial^2}{\partial w^2}-w( w'\frac{\partial^2}{\partial
     w' \partial
     w}+w''\frac{\partial^2}{\partial
     w'' \partial
     w})+\left[2\beta+(-\alpha+\beta+1)w\right]\frac{\partial}{\partial
     w}-\beta' w' \frac{\partial}{\partial w'}-\beta'' w'' \frac{\partial}{\partial w''}+\alpha\beta]\varphi(w, w', w'')=0,
   \\& [w'(1-w')\frac{\partial^2}{\partial w'^2}-w'( w\frac{\partial^2}{\partial
     w \partial
     w'}+w''\frac{\partial^2}{\partial
     w'' \partial
     w'})+\left[2\beta'+(-\alpha+\beta'+1)w'\right]\frac{\partial}{\partial
     w'}-\beta w \frac{\partial}{\partial w}-\beta'' w'' \frac{\partial}{\partial w''}+\alpha\beta']\varphi(w, w', w'')=0,
 \\&[w''(1-w'')\frac{\partial^2}{\partial w''^2}-w''( w'\frac{\partial^2}{\partial
     w'' \partial
     w'}+w\frac{\partial^2}{\partial
     w'' \partial
     w})+\left[2\beta''+(-\alpha+\beta''+1)w''\right]\frac{\partial}{\partial
     w''}-\beta w\frac{\partial}{\partial w}-\beta' w' \frac{\partial}{\partial w'}+\alpha\beta'']\varphi(w, w', w'')=0,
   \end{aligned}$\\
This is an $F_A^{(3)}$  three variable Laurichella hypergeometric system and for $2\beta\neq 1$ and
$2\beta'\neq 1$  $2\beta''\neq 1$  the system has six independent solutions of the form  $[4], p.150-151$:
\begin{itemize}
\item $F_A^{(3)}(-\alpha,\beta,\beta', \beta'', 2\beta,2\beta', 2\beta'', w, w', w'')$,
\item $w^{1-2\beta}F_A^{(3)}(-\alpha+1-2\beta, 1-\beta, \beta', \beta'', 2-2\beta, 2\beta', 2\beta'', w, w', w'')$,
\item $w'^{1-2\beta'}F_A^{(3)}(-\alpha+1-2\beta',\beta, 1-\beta', \beta'', 2\beta, 1-2\beta', 2\beta'', w, w', w'')$,
\item $w''^{1-2\beta''}F_A^{(3)}(-\alpha+1-2\beta'',\beta,\beta', 1-\beta'', 2\beta,2\beta', 2-2\beta'', w, w', w'')$,
\item $w^{1-2\beta}w'^{1-2\beta'}F_A^{(3)}(-\alpha+2-2\beta-2\beta',1-\beta, 1-\beta', \beta'', 2-2\beta, 2-2\beta', 2\beta'', w, w', w'')$,
\item $w^{1-2\beta}w''^{1-2\beta''}F_A^{(3)}(-\alpha+2-2\beta-2\beta'',1-\beta,\beta', 1-\beta'', 2-2\beta, 2\beta',2- 2\beta'', w, w', w'')$,
\item $w'^{1-2\beta}w''^{1-2\beta}F_A^{(3)}(-\alpha+2-2\beta'-2\beta'',\beta,1-\beta', 1-\beta'', 2\beta, 2-2\beta', 2-2\beta'', w, w', w'')$,
\item $ w^{1-2\beta}w'^{1-2\beta}w''^{1-2\beta}F_A^{(3)}(-\alpha +3-2\beta-2\beta'-2\beta'',1-\beta, 1-\beta', 1-\beta'', 2-2\beta, 2-2\beta', 2-2\beta'', w, w', w'')$\\
\end{itemize}
And the proof of the theorem $1.1$ is finished.\\
\section{Cauchy problem for the wave equation with the three-inverse square potential on $R_+^3$}
{\bf Proof of the Theorem B}\\
{\bf Lemma 3.1}{ \it
 Let $F_A^{(3)}$ be the  Appell hypergeometric function with $(h, k, l)\in R^3$ and $a\in R^\ast$
then we have:\\
i)\\
$
  \frac{d}{da}\left[a^{\alpha}F_A^{(3)}\left(-\alpha, \beta, \beta', \beta'', 2\beta,
2\beta', 2\beta'', h/a, k/a, l/a\right)\right]=-\alpha a^{\alpha-1}\times$
$$F_A^{(3)}\left(-\alpha+1, \beta,\beta', \beta'', 2\beta, 2\beta', 2\beta'', h/a,
    k/a, l/a\right)\eqno(3.1)$$
ii)\\
$
a^{\alpha}\Gamma(-\alpha)F_A^{(3)}\left(-\alpha,\beta,\beta', \beta'', 2\beta,
  2\beta',  2\beta'', h/a, k/a, l/a\right)=\Gamma(\beta+1/2)\Gamma(\beta'+1/2)\Gamma(\beta''+1/2)\times$
$$\int_0^\infty e^{-(a-h/2- k/2-l/2)t}t^{-\alpha-1}\left(\frac{th}{4}\right)^{1/2-\beta}\left(\frac{tk}{4}\right)^{1/2-\beta'}\left(\frac{tk}{4}\right)^{1/2-\beta''}\times$$$$I_{\beta-1/2}(th/2)
 I_{\beta'-1/2}(tk/2) I_{\beta''-1/2}(t l/2)dt,
\eqno(3.2)$$
iii) \\
  $
    a^{\alpha}\Gamma(-\alpha)F_A^{(3)}\left(-\alpha,\beta,\beta', \beta''
      \gamma,\gamma', \gamma'', h/a, k/a, l/a\right)\sim$
    $$\Gamma(-\beta-\beta'-\beta''-\alpha)\frac{\Gamma(2\beta)}{\Gamma(\beta)}
    \frac{\Gamma(2\beta')}{\Gamma(\beta')}\frac{\Gamma(2\beta'')}{\Gamma(\beta'')}h^{-\beta}k^{-\beta'}l^{-\beta''}\left(a-h-k-l\right)_+^{\alpha+\beta+\beta'+\beta''}
  \eqno(3.3)$$
as $a\rightarrow 0$.
}.\\
{\bf Proof}\\
i)
is a consequence of the formulas\\
$
a^\alpha
F_A^{(3)}\left(-\alpha,\beta,\beta', \beta'', \gamma, \gamma', \gamma'', h/a, k/a, l/a\right)=$$$\sum_{m, n, p\geq
  0}\frac{(-\alpha)_{m+n+p}(\beta)_m(\beta')_n(\beta'')_p} {(\gamma)_m m!(\gamma')_n n!(\gamma'')_p
  p!}h^m k^n l^pa^{-m-n-p+\alpha}
\eqno(3.4)$$
$
\frac{d}{da}[a^\alpha F_A^{(3)}]\left(-\alpha, \beta, \beta', \beta'', \gamma,
  \gamma', \gamma'', h/a, k/a, l/h\right)]=-\alpha a^{\alpha-1}\times$$$\sum_{m, n, p\geq
  0}\frac{(-\alpha+1)_{m+n+p}(\beta)_m(\beta')_n(\beta'')_p}{(\gamma)_m(\gamma')_n(\gamma'')_p m!
  n!p!}h^mk^n l^pa^{-m-n-p}
\eqno(3.5)$$
$
\frac{d}{da}[a^\alpha F_A^{(3)}]\left(-\alpha, \beta, \beta', \gamma,
\gamma', \gamma'', h/a, k/a, l/a\right)]=-\alpha a^{\alpha-1} \times$$$ F_A\left(-\alpha+1, \beta, \beta', \beta'',
\gamma, \gamma', \gamma'', h/a, k/a, l/a\right)
\eqno(3.6)$$
To prove ii) we use the formulas($[8]$, p.$237$)
$$I_\nu(z)=\frac{2z)^{\nu}e^{-z}}{\sqrt{\pi}\Gamma(\nu+1/2)}\int_0^1e^{2zt}\left[t(1-t)\right]^{\nu-1/2} dt\eqno(3.7)$$
where $\alpha\in C,\beta, \gamma,z\in R^n$
$[1], p.115$
$$F_A^{(3)}\left(\alpha,\beta,\gamma, z\right)=c\int_0^1\int_0^1\int_0^1\prod_{j=1}^{j=3}(1-u_j)^{\gamma_j-\beta_j-1}
u_j^{\beta_j-1}\left(1-\sum_{j=1}^{j=3}u_jz_j\right)^{-\alpha}du_j\eqno(3.8)$$
where $$c=\frac{\prod_{j=1}^{j=3}\Gamma(\gamma_j)}{\prod_{j=1}^{j=3}\Gamma(\beta_j)\Gamma(\gamma_j-\beta_j)}\eqno(3.9)$$
iii) The Proof of iii) uses ii) and the formula ($[8]$, p.$240$)
$$I_{\beta_j-1/2}(x)=\frac{\Gamma(2\beta_j)}{\Gamma(\beta_j)\Gamma(\beta_j+1/2)}2^{-2\beta_j+1/2}x^{-1/2}
\left(1+O(|x|^{-1})\right)\eqno(3.10)$$
To finish the proof of the theorem 1.2, we prove the limit conditions in $(b)$:
from iii) of the Lemma 2.1 and the Legendre duplication formula 
for the $\Gamma$-Euler
function $[6], p.3$
$$
  2^{2z-1}\Gamma(z)\Gamma(z+1/2)=\sqrt{\pi}\Gamma(2z).
\eqno(3.11)$$
we
have For $t \longrightarrow 0$ and $p, p'\in {\mathbb{R}_+^\ast}^3$
$$
  W^{\mathbb{R}^3}_{(b, b', b'')}(t, p, p')\sim c_3 1_{\{|p-p'|<t\}}.
\eqno(3.12)$$
The polar coordinates $p'=p+r\omega$ for $t
\longrightarrow 0$
$$
  u(t,p)\sim c_3\frac{\partial }{t \partial t}\int
  _{{\mathbb{R}^{\ast}}^3}1_{\{|p-p'|<t\}} f(p')dp'
\eqno(3.13)$$
$$
  u(t, p)\sim c_3\frac{\partial }{t \partial t}\int_{0}^t f_p^\#(r)r^2dr \eqno(3.14)$$

$$
  f_p^\#(r)=\frac{1}{6}\int_{S^2 }f(p+r\omega)d\omega
\eqno(3.15)$$
set $r=st$ in the expression above 
 and we see that the limit condition $(b)$ is satisfied and
the proof of the theorem B is finished.

\begin{flushleft}Department of Mathematics\\
College of Arts and Sciences at Al Qurayat   \\
 Al-Jouf University Kingdom of Saudi Arabia  \\
and\\
D\'epartement de Math\'ematiques et Informatique\\
 Facult\'e des Sciences et Techniques\\
Universit\'e de  Nouakchott Al-Aasriya\\
Nouakchott-Mauritanie\\
Email: mohamedvall.ouldmoustapha230@gmail.com
\end{flushleft} 

\end{document}